\documentclass[11pt,a4paper]{article}
\usepackage{amssymb, amsmath, fullpage, mathrsfs, framed, enumerate}
\usepackage{color}
\usepackage{float}
\usepackage{hhline, geometry}
\usepackage{latexsym}
\usepackage{amsmath}
\usepackage{amssymb}
\usepackage{amsthm}
\usepackage{hyperref}
\usepackage{graphicx}
\usepackage{color}
\usepackage{multirow,bigdelim}
\usepackage{natbib}
\usepackage{doi}

\newtheorem{theorem}{Theorem}
\newtheorem{theorem*}{Theorem}
\newtheorem{corollary}[theorem]{Corollary}
\newtheorem{lemma}[theorem]{Lemma}

\theoremstyle{definition}

\newcommand{\GL}{GL}
\newcommand{\F}{\mathbb{F}}
\newcommand{\Z}{\mathbb{Z}}
\newcommand{\Q}{\mathbb{Q}}
\newcommand{\C}{\mathbb{C}}
\newcommand{\N}{\mathbb{N}}

\newcommand{\reg}{\mathrm{reg}}
\newcommand{\regsp}[1]{{#1}_{\reg}}

\newcommand{\rk}{\mathrm{rk}}

\newcommand{\im}{\mathrm{im}}

\newcommand{\codim}{\mathrm{codim}}

\newcommand{\linspan}{\mathrm{span}}
\newcommand{\nc}{\newcommand}
\nc{\rnc}{\renewcommand}
\nc{\beg}{\begin{equation}}
\nc{\eeq}{{\end{equation}}}
\nc{\beqa}{\begin{eqnarray}}
\nc{\eeqa}{\end{eqnarray}}
\nc{\lbar}[1]{\overline{#1}}
\nc{\bra}[1]{\langle#1|}
\nc{\ket}[1]{|#1\rangle}
\nc{\ketbra}[2]{|#1\rangle\!\langle#2|}
\nc{\braket}[2]{\langle#1|#2\rangle}

\nc{\proj}[1]{| #1\rangle\!\langle #1 |}
\nc{\avg}[1]{\langle#1\rangle}
\nc{\Rank}{\operatorname{Rank}}
\nc{\smfrac}[2]{\mbox{$\frac{#1}{#2}$}}
\nc{\tr}{\operatorname{Tr}}
\nc{\ox}{\otimes}
\nc{\dg}{\dagger}
\nc{\dn}{\downarrow}
\nc{\cA}{{\cal A}}
\nc{\cB}{{\cal B}}
\nc{\cC}{{\cal C}}
\nc{\cD}{{\cal D}}
\nc{\cE}{{\cal E}}
\nc{\cF}{{\cal F}}
\nc{\cG}{{\cal G}}
\nc{\cH}{{\cal H}}
\nc{\cI}{{\cal I}}
\nc{\cJ}{{\cal J}}
\nc{\cK}{{\cal K}}
\nc{\cL}{{\cal L}}
\nc{\cM}{{\cal M}}
\nc{\cN}{{\cal N}}
\nc{\cO}{{\cal O}}
\nc{\cP}{{\cal P}}
\nc{\cQ}{{\cal Q}}
\nc{\cR}{{\cal R}}
\nc{\cS}{{\cal S}}
\nc{\cT}{{\cal T}}
\nc{\cX}{{\cal X}}
\nc{\cZ}{{\cal Z}}
\nc{\csupp}{{\operatorname{csupp}}}
\nc{\qsupp}{{\operatorname{qsupp}}}
\nc{\var}{{\operatorname{var}}}
\nc{\rar}{\rightarrow}
\nc{\lrar}{\longrightarrow}
\nc{\polylog}{{\operatorname{polylog}}}
\nc{\1}{{\mathds{1}}}
\nc{\wt}{{\operatorname{wt}}}
\nc{\av}[1]{{\left\langle {#1} \right\rangle}}

\nc{\cAreg}{{{\cal A}_{reg}}}
\nc{\RNS}{{\rm RNS}}
\nc{\RND}{{\rm RND}}

\floatstyle{boxed}
\newfloat{algorithm}{thp}{lop}
\floatname{algorithm}{Algorithm}

\begin{document}

\title{On rank-critical matrix spaces}

\author{
Yinan Li \thanks{Centre for Quantum Software and Information, 
 University of Technology Sydney, Australia  ({\tt liyinan9252@gmail.com}).}
\and 
Youming Qiao \thanks{Centre for Quantum Software and Information, 
 University of Technology Sydney, Australia  ({\tt jimmyqiao86@gmail.com}).}
}
\date{\today}
\maketitle

\begin{abstract}
A matrix space of size $m\times n$ is a linear subspace of the linear space of 
$m\times n$ matrices 
over a field $\F$. The rank of a matrix space is defined as the maximal rank over 
matrices in this space. A matrix space $\cA$ is called rank-critical, if any 
matrix space which properly contains it has rank strictly greater than that of $\cA$. 

In this note, we first exhibit a necessary and sufficient condition 
for a matrix space $\cA$ to be rank-critical, when $\F$ is large enough. 
This immediately implies the sufficient condition for a matrix space to 
be rank-critical by Draisma (Bull. Lond. Math. Soc. 38(5):764--776, 2006), 
albeit requiring the field to be slightly larger.

We then study rank-critical spaces in the context of compression and primitive 
matrix spaces.
We first show that every rank-critical 
matrix space can be decomposed into a rank-critical compression matrix space and a 
rank-critical primitive matrix space. 
We then prove, using our necessary and sufficient 
condition, that the block-diagonal direct sum of two 
rank-critical matrix spaces is rank-critical if and only if both matrix spaces are 
primitive, 
when the field is large enough. 

\end{abstract}

\section{Results}

\subsection{A necessary and sufficient condition for a matrix space to be 
rank-critical}
Let $\F$ be a field, and let $M(m\times n, \F)$ be the linear space of $m\times n$ 
matrices over $\F$. A {\em matrix space} $\cA$ is a linear subspace of $M(m\times 
n, \F)$, 
denoted as $\cA\leq M(m\times n, \F)$. For $A\in M(m\times n, \F)$ we denote its 
rank, kernel, and image, by $\rk(A)$, $\ker(A)$, and $\im(A)$, respectively. 
The \emph{rank of a matrix space} $\cA$, 
denoted as  
$\rk(\cA)$, is defined as $\max\{\rk(A) : A\in 
\cA\}$. $\cA$ is \emph{singular}, if $\rk(\cA)<\min\{m, n\}$. $\cA$ is called {\em 
rank-critical}, if for any $\cB\leq M(m\times n, \F)$ with $\cB\supsetneqq \cA$, 
$\rk(\cB)>\rk(\cA)$. Every $(g, h)\in \GL(m, \F)\times \GL(n, \F)$ has a natural 
action on matrix spaces in $M(m\times n, \F)$, by sending $\cA$ to $g\cA h^{-1}$. 
Two matrix spaces are equivalent if they are in the same orbit of this action. 

Our first result is a necessary and sufficient condition for a matrix space 
to be rank-critical. To state it, we introduce some notation. 
For 
$\cA\leq M(m\times n, \F)$,  
$\regsp{\cA}:=\{A\in \cA : \rk(A)=\rk(\cA)\}$. 
For two subspaces $U\leq  \F^m$, $V\leq\F^n$ and $A\in M(m\times n,\F)$, 
$A(V)=\{A(v) : v\in V\}\leq \F^m$ and $A^{-1}(U)=\{v\in \F^n : A(v)\in U\}\leq 
\F^n$.
Note 
that 
the $A^{-1}$ as in $A^{-1}(\cdot)$ does not refer to the inverse of $A$ and 
$A$ is not necessarily invertible. 

The central notion in our condition is the following. Define the {\em 
rank neutral set} of $\cA\leq M(m\times n, \F)$ as
\begin{equation}
\RNS(\cA):=\{B\in M(m\times n, \F) : \forall A\in \regsp{\cA}, \forall k\in \{0, 
1, \dots, 
m\}, B(A^{-1}B)^{k}\ker(A)\subseteq 
\im(A)\}.
\end{equation}
The elements of $\RNS(\cA)$ are called the {\em rank neutral elements} of $\cA$. 
Note that for $(g, 
h)\in\GL(m, \F)\times \GL(n, \F)$, $\RNS(g\cA h^{-1})=g\RNS(\cA)h^{-1}$. 

\begin{theorem}\label{thm:main}
Let $\cA\leq M(m\times n, \F)$ and suppose $|\F|\geq 2\cdot \min(m, n)$. 
Then $\RNS(\cA)\supseteq \cA$, and $\cA$ is rank-critical if and only if 
$\RNS(\cA)=\cA$. Furthermore, given $G\leq \GL(m, \F)\times\GL(n, \F)$ with the 
natural action on matrix spaces, if $\cA$ is stable under $G$, then 
$\RNS(\cA)$ is also stable under $G$.
%
\end{theorem}

We deduce the sufficient condition for a matrix space to be rank-critical by %
\citet{Dra06}, which plays a key role there to prove that the images of 
certain 
Lie algebra representations are rank-critical. The key notion in Draisma's 
condition is the set of {\em rank neural directions} of $\cA\leq M(m\times n,\F)$,
$$
\RND(\cA):=\{B\in M(m\times n, \F) : \forall A\in \regsp{\cA}, B\ker(A)\subseteq 
\im(A)\}.
$$
Clearly, $\RND(\cA)\supseteq \RNS(\cA)$. Furthermore if 
a group action is present as described in Theorem~\ref{thm:main}, then $\RND(\cA)$ 
is also a stable set under the action of $G$. Therefore the 
following result by Draisma follows immediately from Theorem~\ref{thm:main}. 
\begin{corollary}[{\cite[Proposition 3]{Dra06}}]\label{cor:draisma}
Let $\cA\leq M(m\times n, \F)$ and suppose $|\F|\geq 2\cdot \min(n,m)$.
Then $\RND(\cA)\supseteq \cA$, and if $\RND(\cA)=\cA$ then $\cA$ is rank-critical.
Furthermore, given $G\leq \GL(m, \F)\times\GL(n, \F)$ with the 
natural action on matrix spaces, if $\cA$ is stable under under $G$, then 
$\RND(\cA)$ is also stable under $G$. 
\end{corollary}
We note the following differences between Corollary~\ref{cor:draisma} and 
\citep[Prop. 3]{Dra06}, though such differences are mostly superficial. 
On one hand, Corollary~\ref{cor:draisma} requires the field to be slightly larger 
than needed in \citep[Prop. 3]{Dra06}: there it only requires $|\F|>\rk(\cA)$. On 
the other 
hand, 
Corollary~\ref{cor:draisma} deals with matrix spaces that are not necessarily 
square, and handles a more general group action.

In \citep{Dra06}, Draisma asked the question to investigate the 
``discrepancy between rank-criticality and $\cA=\RND(\cA)$.'' Our result may be 
used as a guide to answer this question: it is now enough to investigate the 
discrepancy between $\RNS(\cA)$ and $\RND(\cA)$. 
Of course, since the condition in the definition of $\RND(\cA)$ is linear, in 
practice it is usually easier to work with $\RND(\cA)$. In fact, we are not aware 
of an explicit example of rank-critical spaces for which the $\RND(\cA)=\cA$ 
fails. 

\subsection{Rank-critical matrix spaces and primitive matrix spaces}

\citet{AL81} introduced the notion of primitive matrix 
spaces. Recall that a matrix space of size $m\times n$ is non-degenerate, if 
$\cap_{A\in 
\cA}\ker(A)=\{0\}$ and $\linspan\{\cup_{A\in 
\cA}\im(A)\}=\F^m$. A matrix space $\cA\leq M(m\times n,\F)$ is  
\begin{itemize}
\item row-primitive, if $\cap_{A\in \regsp{\cA}}\im(A)=0$;
\item column-primitive, if $\linspan\{\cup_{A\in \regsp{\cA}}\ker(A)\}=\F^n$;
\item pre-primitive, if $\cA$ is row-primitive and column-primitive; 
\item primitive, if $\cA$ is non-degenerate, row-primitive, and column-primitive.
\end{itemize}
Note that the zero space in $M(m\times n, \F)$ is also a pre-primitive matrix space.

Another interesting family of matrix spaces is the following. Given $U\leq \F^m$ 
and $V\leq \F^n$, let $p=\dim(U)$, $q=\codim(V)=n-\dim(V)$, and $\cC_{U\leftarrow 
V}^{m, n}
=\{A\in M(m\times n, \F) : A(V)\leq U\}$. 
When $p+q<\min(m, n)$, 
$\rk(\cC_{U\leftarrow V}^{m, n})=\dim(U)+\codim(V)$. 
We call $\cC_{U\leftarrow 
V}^{m, 
n}$ a maximal compression matrix 
space of parameter $(p, q, m, 
n)$. A matrix space is called a compression matrix  space, if it is a subspace of a 
maximal compression matrix space of parameter $(p, q, m, n)$, and its rank is $p+q$. The 
standard maximal compression matrix space of parameter $(p, q, m, n)$ is 
$\cC_{U'\leftarrow V'}^{m, n}$ where $U'$ is spanned by the first $p$ standard 
basis vector of $\F^m$ and $V'$ is spanned by the last $(n-q)$ standard basis 
vector of $\F^n$. We shall denote it $\cC_{p, q}^{m, n}$ for short. Clearly,
$\cC_{p, q}^{m, n}=\{A\in M(m\times n, \F) : \forall p<i\leq m, q<j\leq n, A(i, 
j)=0\}$, where $A(i,j)$ denotes the $(i,j)$th entry of $A$. The standard 
complement of $\cC_{p, q}^{m, n}$ is 
$\overline{\cC_{p, q}^{m, 
n}}:=\{A\in M(m\times n, \F) : \forall 1\leq i\leq p, 1\leq k\leq n, A(i, 
k)=0, \text{ and } \forall 1\leq k\leq m, 1\leq j\leq q, A(k, j)=0\}$. Note 
that $\cC_{0,0}^{m,n}$ is the 
zero matrix space.

%

The following structural result regarding matrix 
spaces of rank bounded from above was first observed by 
\cite{AL81}. 
\begin{theorem}[{\cite[Theorem 1]{AL81}}]\label{thm:decomposition}
Given a singular matrix space $\cA\leq M(m\times n,\F)$, there exist integers 
$p,q\geq 0$ satisfying $p+q< \min(m,n)$, and a primitive matrix space 
$\cP\leq 
M(r\times s,\F)$, $0\leq r\leq m-p$ and $0\leq s\leq n-q$, such that 
$\rk(\cA)=\rk(\cP)+p+q$. Moreover, $\cA$ is equivalent to a matrix space $\cB$ 
in which 
each matrix is of the form
\begin{equation}\label{canonical}
\left[
\begin{array}{c|cc}
p\times q &  & \\ \hline
 & P & 0 \\
& 0 & 0 
\end{array}
\right],
\end{equation}
where $P\in\cP$. 
\end{theorem}
Some remarks are due for this theorem. Firstly, the parameters $p$, $q$, and 
$\cP$ 
are not unique for a given $\cA$. Secondly, the existence of some $p$, $q$, 
and $\cP$ is easy to prove by induction. The main contribution of \cite{AL81} was 
to obtain strong 
restrictions on the size of a primitive matrix space in terms of its rank. Thirdly, when 
$p=q=0$, then $\cA$ is pre-primitive. On the other hand, if $r=s=0$, then $\cA$ is 
a compression matrix space. 


We then study rank-critical matrix spaces 
in the context of 
Theorem~\ref{thm:decomposition}. We first observe that a compression matrix space is 
rank-critical, if and only if it is a maximal compression matrix space (see e.g. 
\citep[Example 10]{Dra06}). In general, for any rank-critical matrix space we have 
the following. 

\begin{theorem}\label{thm:decomposition_critical}
Let $\cA\leq M(m\times n, \F)$ be a matrix space and let $\cB$, $\cP$ be matrix 
spaces as in Theorem~\ref{thm:decomposition}. Let $\cB_c$ be the projection of 
$\cB$ to $\cC_{p, q}^{m,n}$ along $\overline{\cC_{p, q}^{m,n}}$, and $\cB_p$ the 
projection of $\cB$ to $\overline{\cC_{p, q}^{m,n}}$ along $\cC_{p, q}^{m,n}$.
Then $\cA$ is 
rank-critical, if and only if the following hold: (1) $\cB_c=\cC_{p, q}^{m,n}$,
(2) $\cP$ is rank-critical, 
and (3) $\cB= \cB_p\oplus \cB_c$, where $\oplus$ denotes the direct sum of two 
subspaces in $M(m\times n, \F)$. 
\end{theorem}


When $\F$ is large enough, in Theorem~\ref{thm:decomposition_critical}, we may 
replace ``rank-critical'' with the condition $\RNS(*)=*$. It is then 
interesting to consider an analogous statement with $\RND$ instead of $\RNS$.
\begin{theorem}\label{thm:decomposition_rnd}
Suppose $|\F|\geq 2\cdot \min(m, n)$, 
and let $\cA\leq M(m\times n, \F)$ be a 
matrix 
space and let $\cB$, $\cP$ be matrix 
spaces as in Theorem~\ref{thm:decomposition}. Let $\cB_c$ be the projection of 
$\cB$ to $\cC_{p, q}^{m,n}$ along $\overline{\cC_{p, q}^{m,n}}$, and $\cB_p$ the 
projection of $\cB$ to $\overline{\cC_{p, q}^{m,n}}$ along $\cC_{p, q}^{m,n}$.
Then $\RND(\cA)=\cA$, if and only if the following hold: (1) $\cB_c=\cC_{p, q}^{m,n}$,
(2) $\RND(\cP)=\cP$, 
and (3) $\cB= \cB_p\oplus \cB_c$, where $\oplus$ denotes the direct sum of two 
subspaces in $M(m\times n, \F)$. 
\end{theorem}
Theorem~\ref{thm:decomposition_rnd} confirms the common wisdom that to find a 
rank-critical matrix space $\cA$ with $\cA\neq \RND(\cA)$, it is enough to focus 
on primitive matrix spaces. 

Finally, we apply the necessary and sufficient condition from 
Theorem~\ref{thm:main} to prove the following 
result
concerning direct sums of rank-critical matrix spaces. Given two 
matrix spaces $\cA_1\leq M(m_1\times n_1,\F)$ and 
$\cA_2\leq M(m_2\times n_2,\F)$, the (block-diagonal) direct sum of $\cA_1$ and 
$\cA_2$ is a matrix space in $M((m_1+m_2)\times (n_1+n_2), \F)$, defined as 
$\{\left[
\begin{array}{cc}
A_1 &  0 \\
0 & A_2
\end{array}
\right] \in M((m_1+m_2)\times (n_1+n_2), \F): A_1\in \cA_1, A_2\in \cA_2 \}$. By 
abuse of notation we also denote this by $\cA_1\oplus \cA_2$. 

\begin{theorem}\label{thm:sum}
Suppose we are given two rank-critical matrix spaces $\cA_1\leq M(m_1\times 
n_1,\F)$ and 
$\cA_2\leq M(m_2\times n_2,\F)$, and suppose $|\F|\geq 2\min(m_1+m_2, n_1+n_2)$. $\cA_1\oplus\cA_2$ is rank-critical if and only 
if $\cA_1$ and $\cA_2$ are primitive.
\end{theorem}

\section{Proofs}
\subsection{On Theorem~\ref{thm:main}}
\subsubsection{The Wong sequences, and some digression}\label{subsubsec:wong}

Our condition is achieved via
a perspective that is different from Draisma's as in \citep{Dra06}. Draisma 
arrived at the sufficient condition
$\RND(\cA)=\cA$ from a geometric perspective, by considering tangent spaces at 
regular points in a linear subspace contained in an affine variety. On the other 
hand, our condition, $\RNS(\cA)=\cA$, was obtained from an algorithmic 
perspective. 
We now introduce some 
previous results from 
\cite{IKQS15} that support the proof of Theorem~\ref{thm:main}, together with some 
background information. Some of the material here is 
more general than strictly needed to prove Theorem~\ref{thm:main}, as we want to 
take this chance to advocate a connection between the geometry of matrix spaces 
and a key algorithmic problem in computational complexity theory. 

A central problem in computational complexity theory is the {\em symbolic 
determinant 
identity testing} (SDIT) problem, which asks to decide whether a matrix space, 
given by a linear basis, contains a full-rank matrix. When the underlying field is 
large enough, SDIT admits a randomized efficient algorithm \cite{Lov79}. The goal 
then is to devise a deterministic efficient algorithm, as this 
implies an arithmetic circuit lower bound that is believed to be beyond current 
techniques \citep{CIKK15}. 

In fact, for the purpose of \citep{CIKK15}, it is enough to exhibit a 
polynomial-size witness for the singularity of a matrix space. This problem is 
wide open, while some helpful structures are known. One such structure is the 
following. For $\cA\leq M(m\times n, \F)$, and $V\leq \F^n$, it is easy to verify 
that $\rk(\cA)\leq n-(\dim(V)-\dim(\cA(V)))=\codim(V)+\dim(\cA(V))$. So 
$\rk(\cA)\leq n-\max\{\dim(V)-\dim(\cA(V)) : V\leq \F^n\}$.
\citet{Lovasz} observed that if $\cA$ has a basis consisting of rank-$1$ matrices, 
then this upper bound can be achieved at some $V\leq \F^n$. This follows from the 
matroid intersection theorem for linear matroids \citep{Edmonds70}.

Furthermore, for $\cA\leq M(m\times n, \F)$ and $s\in \Z^+$, we call $V\leq \F^n$ an 
$s$-shrunk 
subspace of $\cA$, if $\dim(V)-\dim(\cA(V))\geq s$. It is then an interesting 
question 
to decide whether a given 
matrix space possesses an $s$-shrunk subspace for a given $s\in \Z^+$. Recently, 
deterministic 
polynomial-time 
algorithms were devised in \cite{GGOW} over $\Q$, and \cite{IQSnote,IQS1} 
over any field. 
The key algorithmic technique in \cite{IQSnote,IQS1} is the (second 
generalized) Wong sequences, first used in \cite{FR04} and then rediscovered in 
\cite{IKQS15}. 
They can be viewed as a linear algebraic analogue of the augmenting 
paths, which were developed to solve the perfect matching problem on bipartite 
graphs. 
Given $A\in \cA\leq M(m\times n, \F)$, the Wong sequence of $(A, 
\cA)$ is the following sequence of subspaces of $\F^m$: $W_0=\{0\}, W_1=\cA (A^{-1} 
(W_0)), W_2=\cA (A^{-1} (W_1)), \dots, W_{i+1}=\cA (A^{-1} (W_i)), \dots$. It is 
known 
that for some $\ell\in \{0, 1, \dots, m\}$, 
$W_0<W_1<\dots 
<W_\ell=W_{\ell+1}=\dots$ \cite[Prop. 7]{IKQS15}, and $\cA$ has a 
$\dim(\ker(A))$-shrunk subspace if and only 
if $W_\ell\subseteq \im(A)$ \cite[Lemma 9]{IKQS15}.

\subsubsection{Proof of Theorem~\ref{thm:main}} 

We now turn to prove Theorem~\ref{thm:main}. The reader probably has noticed 
the 
similarity between the formulation of the rank neutral set, and Wong sequences 
introduced above. One more ingredient is to relate Wong sequences to 2-dimensional 
matrix spaces. For $\cA\leq M(m\times n, \F)$ of dimension $2$ with $|\F|> 
\min(m, n)$, it is known that 
$\rk(\cA)=n-\max\{\dim(V)-\dim(\cA(V)) : V\leq \F^n\}$ (see e.g. \cite{AS78}).
Combining with the Wong sequences, \citet{IKQS15} showed the following: 
\begin{lemma}[{\cite[Lemma 12]{IKQS15}}]\label{lem:2dim}
Suppose we are given $A\in \cA=\linspan\{A, B\}\leq M(m\times n, \F)$, and 
$|\F|>\min(m,n)$. Then $A$ is of 
maximal rank in $\cA$, if 
and only if for $i\in \{0, 1, \dots, m\}$, $B(A^{-1}B)^{i}\ker(A)\leq \im(A)$. 
\end{lemma}
Given Lemma~\ref{lem:2dim} it is easy to prove Theorem~\ref{thm:main}.

\paragraph{Theorem~\ref{thm:main}, restated} 
Let $\cA\leq M(m\times n, \F)$ and suppose $|\F|\geq 2\cdot \min(m, n)$. 
Then $\RNS(\cA)\supseteq \cA$, and $\cA$ is rank-critical if and only if 
$\RNS(\cA)=\cA$. Furthermore, given $G\leq \GL(m, \F)\times\GL(n, \F)$ with the 
natural action on matrix spaces, if $\cA$ is stable under $G$, then 
$\RNS(\cA)$ is also stable under $G$.

\begin{proof} 
To start with, note that Lemma~\ref{lem:2dim} immediately implies that 
$\RNS(\cA)\supseteq \cA$. 

We first show that $\RNS(\cA)=\cA$ implies that $\cA$ is rank-critical. By 
contradiction suppose there exists a matrix $B\not\in \cA$ s.t. 
$\rk(\linspan\{\cA, B\})=\rk(\cA)$. Then for any $A\in \regsp{\cA}$, $\rk(\linspan\{A, 
B\})=\rk(A)$. 
Lemma~\ref{lem:2dim} tells us that $B\in \RNS(\cA)$, so $\cA$ is a proper subset 
of $\RNS(\cA)$, a contradiction. 

We then prove that if $\cA$ is rank-critical then $\RNS(\cA)=\cA$. Suppose not, 
then there exists $B\in \RNS(\cA)\setminus \cA$. Let 
$r=\rk(\cA)$, and $\cA'=\linspan\{B, \cA\}$. Note that $r<\min(m, n)$. 
Because $\cA$ is rank-critical, $\rk(\cA')>r$, so there exists $A\in \cA$ 
s.t. $\rk(B+A)>r$. Since $B\in \RNS(\cA)$, by Lemma~\ref{lem:2dim} $A$ 
cannot be from $\regsp{\cA}$. Take any $A'\in \regsp{\cA}$, and consider $B+A+xA'$, where 
$x$ is a formal variable.
As $\rk(B+A)> r$, 
for all but at most $(r+1)$ $\lambda\in \F$, 
$\rk(B+A+\lambda A')=\rk(B+A)>r$. As $\rk(A')=r$, for all but at most $r$ $\mu\in \F$, 
$\rk(A+\mu A')=r$. Since $|\F|\geq 2\min(m, n)>2r+1$, there exists some 
$\nu \in\F$, 
such that $\rk(B+A+\nu A')> r$ and $\rk(A+\nu A')=r$. In this case, $A+\nu 
A'\in\regsp{\cA}$, so by Lemma~\ref{lem:2dim} again, this suggests that $B\not\in 
\RNS(\cA)$, a contradiction. 

To see that the statement regarding the group action holds, recall that $\RNS(g\cA 
h^{-1})=g\RNS(\cA) h^{-1}$. 
\end{proof}

\subsection{Proof of Theorem~\ref{thm:decomposition_critical}}
\paragraph{Theorem~\ref{thm:decomposition_critical}, restated} Let $\cA\leq 
M(m\times n, \F)$ be a matrix space and let $\cB$, $\cP$ be matrix 
spaces as in Theorem~\ref{thm:decomposition}. Let $\cB_c$ be the projection of 
$\cB$ to $\cC_{p, q}^{m,n}$ along $\overline{\cC_{p, q}^{m,n}}$, and $\cB_p$ the 
projection of $\cB$ to $\overline{\cC_{p, q}^{m,n}}$ along $\cC_{p, q}^{m,n}$.
Then $\cA$ is 
rank-critical, if and only if the following hold: (1) $\cB_c=\cC_{p, q}^{m,n}$,
(2) $\cP$ is rank-critical, 
and (3) $\cB=\cB_p\oplus \cB_c$, where $\oplus$ denotes the direct sum of two 
subspaces in $M(m\times n, \F)$.

\begin{proof}
As $\cA$ is rank-critical if and only if $\cB$ is rank-critical, we focus on $\cB$ 
in the following. 

We first examine the necessary direction. Recall that from 
Theorem~\ref{thm:decomposition}, we 
have 
$p, q\in \N$ satisfying $p+q<\min(m,n)$, and a primitive matrix space 
$\cP\leq M(r\times s,\F)$ where $r\leq m-p$ and $s\leq n-q$, such that (1) 
$\rk(\cB)=p+q+\rk(\cP)$, and (2) every $B\in \cB$ is of the form 
\begin{equation}\label{eq:compression}
B=\left[ {\begin{array}{*{20}{c}}
{\begin{array}{*{20}{c}}
\ldelim\{{3}{6pt}[p]\\&\\&\\
\end{array}\overbrace {\begin{array}{*{20}{c}}
{\rm{*}}& \cdots &{\rm{*}}\\
 \vdots &{}& \vdots \\
{\rm{*}}& \cdots &{\rm{*}}
\end{array}}^{\rm{q}}}&{\begin{array}{*{20}{c}}
*& \cdots & \cdots &*\\
 \vdots &{}&{}& \vdots \\
*& \cdots & \cdots &*
\end{array}}\\
{\begin{array}{*{20}{c}}
{}
\end{array}\begin{array}{*{20}{c}}
& &*& \cdots &*\\
 & &\vdots &{}& \vdots \\
& & \vdots &{}& \vdots \\
& &*& \cdots &*
\end{array}}&{\begin{matrix}	P_{r\times s}& &0 \\ & & \\ 0& &0\end{matrix}}
\end{array}} \right]_{m\times n},
\end{equation}
where $P\in\cP$. 

Now $\cB$ is rank-critical. We first show that $\cC_{p, q}^{m,n}\leq \cB$, which 
will then establish (1) and (3). Take any $C\in \cC_{p, q}^{m,n}$, and let 
$\cB'=\linspan(C, \cB)$. Any $B'\in \cB$ is also of the form as in 
Equation~\ref{eq:compression}, as $C$ only adds to the $*$ entries. But this gives 
that $\rk(B')\leq p+q+\rk(P)\leq p+q+\rk(\cP)=\rk(\cB)$. Then by the rank 
criticality 
of $\cB$, $C\in \cB$.

%
We then turn to (2) $\cP$ is rank-critical. Suppose $\cP$ is not, then there 
exists some $P'\not\in\cP$ 
satisfying $\rk(\linspan(\cP,{P'}))=\rk(\cP)$. Then let $C\in 
M(m\times n,\F)$ be
$$C=\left[ {\begin{array}{*{20}{c}}
{\begin{array}{*{20}{c}}
\ldelim\{{3}{6pt}[p]\\&\\&\\
\end{array}\overbrace {\begin{array}{*{20}{c}}
{\rm{0}}& \cdots &{\rm{0}}\\
 \vdots &{}& \vdots \\
{\rm{0}}& \cdots &{\rm{0}}
\end{array}}^{\rm{q}}}&{\begin{array}{*{20}{c}}
0& \cdots & \cdots &0\\
 \vdots &{}&{}& \vdots \\
0& \cdots & \cdots &0
\end{array}}\\
{\begin{array}{*{20}{c}}
{}
\end{array}\begin{array}{*{20}{c}}
& &0& \cdots &0\\
 & &\vdots &{}& \vdots \\
& & \vdots &{}& \vdots \\
& &0& \cdots &0
\end{array}}&{\begin{matrix}	P'& &0 \\ & & \\ 0& &0\end{matrix}}
\end{array}} \right]_{m\times n}.
$$
Clearly, $C\not\in\cB$. Now consider $\cB'=\linspan(C, \cB)$. We then have 
$\rk(\cB')\leq p+q+\rk(\linspan(P', \cP))=p+q+\rk(\cP)=\rk(\cB')$. 
This contradicts the rank-criticality of $\cB$, proving that $\cP$ is 
rank-critical.

To show the sufficiency, our strategy is the following. Let $\cB_{\mathrm{ur}}\leq 
M((p+r)\times (n-q), \F)$ be the matrix space that consists of those submatrices 
of size $(p+r)\times (n-q)$ in the 
upper-right corner of $B\in \cB$. We first prove that $\cB_{\mathrm{ur}}$ is 
rank-critical, 
using only the row primitivity of $\cP$. 
We then show that 
as $\cP$ is 
column-primitive, $\cB_{\mathrm{ur}}$ is also column-primitive. This allows us to 
conclude 
that $\cB$ is rank-critical, by applying the 
column version of the argument which proved the rank criticality 
of 
$\cB_{\mathrm{ur}}$.

We first prove that $\cB_{\mathrm{ur}}$ is rank-critical. To start with, note that 
$\rk(\cB_{\mathrm{ur}})=p+\rk(\cP)$, $p+\rk(\cP)<p+r$ (by $\rk(\cP)<r$), and 
$p+\rk(\cP)<n-q$ (by $p+q+\rk(\cP)<\min\{m, n\}$). So $\cB_{\mathrm{ur}}$ is 
singular. 
Suppose we have $C\in M((p+r)\times (n-q), \F)$, $C\not\in\cB_{\mathrm{ur}}$, such 
that 
$\rk(\cB_{\mathrm{ur}})=\rk(\cB'_{\mathrm{ur}})$ where 
$\cB'_{\mathrm{ur}}=\linspan(C, \cB_{\mathrm{ur}})$. As the 
first $p$ rows are 
free in $\cB_{\mathrm{ur}}$, w.l.o.g. we can assume the first $p$ rows of $C$ are 
$0$. 
Write $C$ as  
$C=\begin{bmatrix}0&0\\C_1&C_2\end{bmatrix}_{(p+r)\times(n-q)}$ where $C_1$ is of size 
$r\times s$. We observe that $C_1\in \cP$. If not, by the rank-criticality 
of $\cP$, we would have $\rk(\cB_{\mathrm{ur}})=\rk(\cP)+p<\rk(\linspan(C_1, 
\cP))+p\leq 
\rk(\cB'_{\mathrm{ur}})$, a contradiction. Therefore we can further assume $C$ to 
be of the 
form 
$C=\begin{bmatrix}0&0\\0&C_2\end{bmatrix}_{(p+r)\times(n-q)}$. Consider the matrix space 
$\cC=\linspan(\begin{bmatrix} 0 & C_2
\end{bmatrix}, \begin{bmatrix} \cP & 0 
\end{bmatrix})\leq M(r\times (n-q), \F)$. As before, since 
$\rk(\cB_{\mathrm{ur}})=\rk(\cB'_{\mathrm{ur}})$, 
it is necessary that $\rk(\cC)=\rk(\cP)$. It follows that every 
column of 
$C_2$ is in $\cap_{P\in \regsp{\cP}} \im(P)$. By the row-primitivity of $\cP$, 
$C_2$ has to be the zero matrix. Therefore the whole $C$ is the zero matrix, 
proving that 
$\cB_{\mathrm{ur}}$ is rank-critical. 

We now prove that $\cB_{\mathrm{ur}}$ is column-primitive. As $\cP$ is 
column-primitive, 
$\begin{bmatrix} \cP & 0 
\end{bmatrix}$ is also column primitive. 
Take any $B\in \regsp{(\cB_{\mathrm{ur}})}$. $B$ is 
of the form $\begin{bmatrix}D_1&D_2\\P&0\end{bmatrix}$, and $\ker(B)=\ker 
(\begin{bmatrix}D_1&D_2\end{bmatrix})\cap \ker(\begin{bmatrix} P & 0 
\end{bmatrix})$. $\ker(B)\neq 0$ since $\cB_{\mathrm{ur}}$ is singular. As the 
first $p$ 
rows are free, by choosing appropriate $\begin{bmatrix}D_1&D_2\end{bmatrix}$ we 
can go through all codimension-$p$ subspaces of $\ker(\begin{bmatrix} P & 
0 
\end{bmatrix})$. Now the column-primitivity of $\cB_{\mathrm{ur}}$ follows from 
that of 
$\begin{bmatrix} \cP & 0 
\end{bmatrix}$.
\end{proof}
\subsection{Proof of Theorem~\ref{thm:decomposition_rnd}}

\paragraph{Theorem~\ref{thm:decomposition_rnd}, restated}
Let $\cA\leq M(m\times n, \F)$ be a matrix space and let $\cB$, $\cP$ be matrix 
spaces as in Theorem~\ref{thm:decomposition}. Let $\cB_c$ be the projection of 
$\cB$ to $\cC_{p, q}^{m,n}$ along $\overline{\cC_{p, q}^{m,n}}$, and $\cB_p$ the 
projection of $\cB$ to $\overline{\cC_{p, q}^{m,n}}$ along $\cC_{p, q}^{m,n}$.
Then $\RND(\cA)=\cA$, if and only if the following hold: (1) $\cB_c=\cC_{p, q}^{m,n}$,
(2) $\RND(\cP)=\cP$, 
and (3) $\cB= \cB_p\oplus \cB_c$, where $\oplus$ denotes the direct sum of two 
subspaces in $M(m\times n, \F)$.
\begin{proof}
To start with, note that when $\cA$ and $\cB$ are equivalent, then $\RND(\cA)=\cA$ 
if and only if $\RND(\cB)=\cB$. The proof strategy is similar to the proof of 
Theorem~\ref{thm:decomposition_critical}, while some changes are required to deal 
with $\RND(\cB)$. 

For the sufficiency direction, 
let $\cB_{\mathrm{ur}}\leq M((p+r)\times (n-q), \F)$ be a matrix space that 
consists of those submatrices of 
size $(p+r)\times (n-q)$ in the 
upper-right corner of $B\in \cB$. We will first show that 
$\RND(\cB_{\mathrm{ur}})=\cB_{\mathrm{ur}}$,
using only the row primitivity of $\cP$. 
Then by the column version of the 
argument, we can conclude that $\RND(\cB)=\cB$, as $\cB_{\mathrm{ur}}$ is also 
column-primitive as shown in the proof of 
Theorem~\ref{thm:decomposition_critical}. 

It remains to prove that $\RND(\cB_{\mathrm{ur}})=\cB_{\mathrm{ur}}$. Take any 
$B\in\regsp{(\cB_{\mathrm{ur}})}$. $B$ is of the 
form $B=\begin{bmatrix}
D_1&D_2\\P&0 \end{bmatrix}_{(p+r)\times(n-q)}$, where $D_1\in 
M(p\times s,\F)$, $D_2\in 
M(p\times(n-q-s),\F)$, and $P\in\cP$. Clearly, 
$\ker(B)=\ker(\begin{bmatrix}D_1&D_2\end{bmatrix})
\cap\ker(\begin{bmatrix}P&0\end{bmatrix})$.
Since $\rk(B)=p+\rk(P)$, $\ker(B)$ is a codimension-$p$ subspace in 
$\ker(\begin{bmatrix}P&0\end{bmatrix})$, which implies that 
$\im(B)=\F^{p}\oplus\im(P)$. 
Now let $C$ be a rank neutral direction 
of $\cB_{\mathrm{ur}}$, and put it in the block form 
$C=\begin{bmatrix}C_{11}&C_{12}\\C_{21}&C_{22}\end{bmatrix}$, where $C_{11}\in 
M(p\times s,\F)$, $C_{12}\in M(p\times(n-q-s),\F)$, $C_{21}\in M(r\times s,\F)$ and $C_{22}\in 
M(r\times(n-q-s),\F)$. By the definition of rank neutral directions, we have
\begin{equation}\label{eq:test}
\begin{split}
\begin{bmatrix} 
C_{11}&C_{12}\end{bmatrix}\big(\ker(\begin{bmatrix}D_1&D_2\end{bmatrix})\cap\ker(\begin{bmatrix}P&0\end{bmatrix})\big)&\leq
 \F^{p};\\
\begin{bmatrix} 
C_{21}&C_{22}\end{bmatrix}\big(\ker(\begin{bmatrix}D_1&D_2\end{bmatrix})\cap\ker(\begin{bmatrix}P&0\end{bmatrix})\big)&\leq
 \im(P)
\end{split}
\end{equation}
for all $D_1$, $D_2$ and $P\in\regsp{\cP}$ satisfying 
$\rk(\begin{bmatrix}D_1&D_2\\P&0\end{bmatrix})=p+\rk(P)$. The first 
constraint in Equation~\ref{eq:test} puts no restriction on  $C_{11}$ 
and $C_{12}$. 
For the second constraint in Equation~\ref{eq:test}, as already argued in the last 
paragraph in the proof of 
Theorem~\ref{thm:decomposition_critical}, since the first $p$ rows are free, by 
choosing appropriate $\begin{bmatrix}D_1&D_2\end{bmatrix}$ we can go over 
all codimension-$p$ 
subspaces of $\ker(\begin{bmatrix}P&0\end{bmatrix})$. 
 This gives that 
$\begin{bmatrix} 
C_{21}&C_{22}\end{bmatrix}\ker(\begin{bmatrix}P&0\end{bmatrix})\leq \im(P)$. 
Then by 
$\ker(\begin{bmatrix}P&0\end{bmatrix})=\ker(P)\oplus \F^{n-q-s}$, 
we have
\begin{equation}\label{eq:test2}
C_{21}\ker(P)+C_{22}\F^{n-q-s}\leq\im(P),
\end{equation}
for all $P\in\regsp{\cP}$, from which we deduce that (a) $C_{21}\in\RND(\cP)=\cP$ 
for any $P\in\regsp{\cP}$, and (b) 
$C_{22}=0$ as $\im(C_{22})\leq 
\cap_{P\in\regsp{\cP}}\im(P)=\{0\}$, where the equality follows from the 
row primitivity of $\cP$. That $\RND(\cB_{\mathrm{ur}})=\cB_{\mathrm{ur}}$ 
then follows.

For the necessary direction, notice that $\RND(\cB)=\cB$ implies 
$\RNS(\cB)=\cB$, thus conditions (1) and (3) hold by Theorem 
\ref{thm:decomposition_critical}. By contradiction, assume that $\RND(\cP)\neq 
\cP$, so there exists 
$P_0\in\RND(\cP)$ but $P_0\not\in\cP$. It is easy to see that 
$P_0'=\begin{bmatrix}0&0\\P_0&0\end{bmatrix}\in M((p+r)\times (n-q), \F)$ is not 
an element of $\cB_{\mathrm{ur}}$ but satisfies Equations~\ref{eq:test} for all $P\in 
\regsp{\cP}$, 
which implies $P'_0\in\RND(\cB_{\mathrm{ur}})$. Consider then the matrix 
$P''_0=\begin{bmatrix}0&0&0\\0&P_0&0\\0&0&0\end{bmatrix}\in 
M(m\times n,\F)$, and by the column version of the argument, we have $P''_0\in\RND(\cB)$ but 
$P''_0\not\in\cB$, arriving at a 
contradiction. 
\end{proof}

\subsection{Proof of Theorem~\ref{thm:sum}}
\paragraph{Theorem~\ref{thm:sum}, restated} 
Suppose we are given two rank-critical matrix spaces $\cA_1\leq M(m_1\times 
n_1,\F)$ and 
$\cA_2\leq M(m_2\times n_2,\F)$, and suppose $|\F|\geq 2\min(m_1+m_2, n_1+n_2)$. 
$\cA_1\oplus\cA_2$ is rank-critical if and only 
if $\cA_1$ and $\cA_2$ are primitive.

We point out that, by the discussion in Section~\ref{subsubsec:wong}, an 
equivalent 
formulation of $\RNS(\cA)$ is 
$$
\RNS(\cA):=\{B\in M(m\times n, \F) : \forall A\in \regsp{\cA}, \forall k\in \N, 
B(A^{-1}B)^{k}\ker(A)\subseteq 
\im(A)\}.
$$

\begin{proof}
To see the necessity, we prove that if $\cA_1$ is not primitive, then $\cA_1\oplus 
\cA_2$ is not rank-critical. If $\cA_1$ is not primitive then 
$\cA_1\oplus \cA_2$ is not primitive. Furthermore by transforming to an equivalent 
space, $\cA_1$ can be arranged to be in the form as the $\cB$ in 
Theorem~\ref{thm:decomposition}, 
with one of $p$ or $q$ being nonzero. W.l.o.g. assume $q>0$. Then the first column 
is also the cause of imprimitivity of $\cA_1\oplus \cA_2$; that is, the first 
standard basis vector is not in $\linspan\{\cup_{A\in \regsp{(\cA_1\oplus 
\cA_2)}}\ker(A)\}$. Now 
by 
Theorem~\ref{thm:decomposition_critical}, for $\cA_1\oplus\cA_2$ to be 
rank-critical, it is necessary that the first column is free, while every $A\in 
\cA_1\oplus 
\cA_2$ would have the first column containing some $0$'s. This proves that 
$\cA_1\oplus 
\cA_2$ is not rank-critical. 


For the sufficiency direction, by Theorem \ref{thm:main}, we turn to prove 
$\cA_1\oplus\cA_2=\RNS(\cA_1\oplus\cA_2)$. That is, for any 
$X\in\RNS(\cA_1\oplus\cA_2)$ satisfying $\forall A\in \regsp{(\cA_1\oplus\cA_2)}, 
\forall k\in\N$,
\begin{equation}\label{condition}
X(A^{-1}X)^k\ker(A)\leq\im(A),
\end{equation}
we need to show $X\in\cA_1\oplus\cA_2$. 
Noticing $\regsp{(\cA_1\oplus\cA_2)}=\regsp{(\cA_1)}\oplus \regsp{(\cA_2)}$, we 
denote a given  
$A\in \regsp{(\cA_1\oplus\cA_2)}$ by $A=A_1\oplus A_2$, where 
$A_1\in \regsp{(\cA_1)}$ and $A_2\in \regsp{(\cA_2)}$. Moreover, we have 
$\im(A)=\im(A_1)\oplus\im(A_2)$ and $\ker(A)=\ker(A_1)\oplus\ker(A_2)$.

Now, let $X=\begin{bmatrix}X_{11}&X_{12}\\X_{21}&X_{22}\end{bmatrix}\in 
\RNS(\cA_1\oplus\cA_2)$, where $X_{ij}\in M(m_i\times n_j,\F)$, $i,j=1,2$. By 
Equation~\ref{condition} with $k=0$, for any $A_1\in \regsp{(\cA_1)}$, $A_2\in 
\regsp{(\cA_2)}$, 
$X_{ij}$'s satisfy:
$$X_{11}\ker(A_1)+X_{12}\ker(A_2)\leq\im(A_1);$$
$$X_{21}\ker(A_1)+X_{22}\ker(A_2)\leq\im(A_2).$$
Therefore, $X_{12}\ker(A_2)\leq\im(A_1)$ and $X_{21}\ker(A_1)\leq\im(A_2)$ hold 
for any $A_1\in \regsp{(\cA_1)}$, $A_2\in \regsp{(\cA_2)}$. So we have
$$X_{12}(\linspan\{\cup_{A_2\in \regsp{(\cA_2)}}\ker(A_2)\})\leq\cap_{A_1\in 
\regsp{(\cA_1)}}\im(A_1);$$
$$X_{21}(\linspan\{\cup_{A_1\in \regsp{(\cA_1)}}\ker(A_1)\})\leq\cap_{A_2\in 
\regsp{(\cA_2)}}\im(A_2).$$
Now by the primitivity of $\cA_1$ and $\cA_2$, we obtain $X_{12}=0$ and $X_{21}=0$.

We then need to show that for $i=1, 2$, $X_{ii}\in\cA_i$. By the assumption 
$\cA_i=\RNS(\cA_i)$, we turn to show that for $i=1, 2$, $X_{ii}\in\RNS(\cA_i)$, 
that is, $\forall k\in\N$, and $i=1, 2$, 
$X_{ii}(A_i^{-1}X_{ii})^k\ker(A_i)\leq\im(A_i)$. This can be seen by an induction 
on $k$, once we notice the following: if $U=U_1\oplus U_2$, $U_i\in \im(A_i)$, 
then $A^{-1}(U)=A_1^{-1}(U_1)\oplus A_2^{-1}(U_2)$. 
\end{proof}

\subsection*{Acknowledgement}

We thank Jan Draisma and G\'abor Ivanyos for helpful discussions though email 
correspondences. Y. Q. 
was supported by the Australian Research Council DECRA DE150100720 during this 
research. 

\section*{Bibliography}

\bibliographystyle{plainnat}
\bibliography{references}

%
%
%

\end{document}